\documentclass[12pt,a4paper,reqno]{amsart}
\usepackage[utf8]{inputenc}
\usepackage{color}
\usepackage{enumerate}
\usepackage[margin=2.5cm]{geometry}
\usepackage{graphicx} 
\usepackage{overpic}
\usepackage{amsaddr}

\newtheorem{lema}{Lemma}[section]
\newtheorem{teo}[lema]{Theorem}
\newtheorem{pro}[lema]{Proposition}
\newtheorem{remark}[lema]{Remark}
\newtheorem{defini}[lema]{Definition}

\DeclareMathOperator{\Res}{Res}

\newcommand{\R}{\ensuremath{\mathbb{R}}}
\newcommand{\C}{\ensuremath{\mathbb{C}}}

\title[Kolmogorov system]{New Lower Bound For the Hilbert Number in low degree Kolmogorov systems}

\author[Y. R. Carvalho] {Yagor Romano Carvalho$^{1}$}
\email{$^1$yagor.carvalho@usp.br}

\author[L. P. C. da Cruz] {Leonardo P. C. da Cruz$^2$}
\address{$^{1,2}$Universidade de São Paulo, São Carlos, Brazil}
\email{$^2$leonardocruz@icmc.usp.br}

\author[L. F. S Golveia] {Luiz F.S. Gouveia$^*$}
\address{$^*$ UNESP, São José do Rio Preto, Brazil }
\email{$^*$fernando.gouveia@unesp.br}



\begin{document}
	
\maketitle

\begin{abstract} 
Our main goal in this paper is to study the number of small-amplitude isolated periodic orbits, so-called limit	cycles, surrounding only one equilibrium point a class of polynomial Kolmogorov systems. We denote by $\mathcal M_{K}(n)$ the maximum number of limit cycles bifurcating from the equilibrium point via a degenerate Hopf bifurcation for a polynomial Kolmogorov vector field of degree $n$. In this work, we obtain another example such that  $ \mathcal M_{K}(3)\geq 6$. In addition, we obtain new lower bounds for $\mathcal M_{K}(n)$ proving that $\mathcal M_{K}(4)\geq 13$ and $\mathcal M_{K}(5)\geq 22$.
\end{abstract}
	
\section{Introduction}

In the years 1925 and 1926 Alfred J. Lotka \cite{Lotka1925} and Vito Volterra \cite{Volterra1927} independently propose a class of polynomial differential systems of degree two in the plane for studying the interactions of two populations (species occupying the same ecological niche), that we call nowadays of Lotka-Volterra systems. Such systems have the following form

\begin{equation}\label{eq:1}
	\begin{array}{lll}
			\dot{x}&=&x\,(a_0 + a_1 x+a_2y),\\
			\dot{y}&=&y\,(b_0 + b_1 x+b_2y),
	\end{array}
\end{equation} 
where $a_i,b_i \in \R, i=0,1,2$. 

When studying the interactions between two populations, it only makes sense to restrict the analysis of the orbit behavior to the first quadrant. In this work, we would like to study bifurcation of limit cycles, see \cite{Albrecht1973,Huang1988,May1972}. A limit cycle is a periodic orbit of a differential system that is isolated in the set of all its periodic orbits.

In 1936 Kolmogorov \cite{Kolmogorov1936} extended the Lotka-Volterra systems to arbitrary degrees and dimensions. Then a planar Kolmogorov differential system is a planar dynamical system of the form
\begin{equation}\label{eq:2}
	\begin{array}{lll}
		\dot{x}&=&x\,X_m(x,y),\\
		\dot{y}&=&y\,Y_m(x,y),
	\end{array}
\end{equation} 
where $X_m$ and $Y_m$ are polynomials of degree $m$. Particularly when $m=1$, we have the Lotka-Volterra systems. As in the Lotka-Volterra system, the equilibrium point is located in the quadrant where the Kolmogorov systems have biological meaning, that is, in the first quadrant. In addition, this class of systems has a wide range of applications, such as chemical reactions \cite{Her1990}, economics \cite{Gandolfo2008,Goodwin1967,SoRic2002} and hydrodynamics \cite{Bus1981}.

The limit cycles have been widely studied in the last century, and the analysis of the existence of these objects is significant in studying the interactions between two species. In general, studying the maximum number of limit cycles of (usually denoted by $\mathcal{H}(n)$ and called by Hilbert number) vector fields in terms of their degrees is an open problem and approached in the second part of the 16th Hilbert Problem, see \cite{Ily2002,Sma1998}. In this work, we are interested in a local version of the 16th Hilbert Problem that consists in providing the maximum number $\mathcal M(n)$ of small-amplitude limit cycles bifurcating from an equilibrium point of center or focus type, and clearly $\mathcal M(n) \leq \mathcal{H}(n)$. The number of small-amplitude limit cycles born by small perturbations from the equilibrium point is usually referred to as the cyclicity of this equilibrium point. In other words, $M(n)$ is an upper bound of the cyclicity of such equilibrium points.  In our context, cyclicity is the maximum number of limit cycles that can bifurcate from the equilibrium point by means of a degenerate Hopf bifurcation. More details on related problems about the cyclicity of homoclinic and heteroclinic connections or period annulus can be found in \cite{TorreGouv2021,Roussarie1998}.

It is well known that linear systems do not have limit cycles, then $\mathcal{H}(n)=0$. For $n=2,$ the problem of estimating $\mathcal{H}(2)$ has been studied intensively during the last century. The lower bounds for $\mathcal{H}(2)$ can be given by providing concrete examples of polynomial differential systems of degree 2. Up to now, the best result was given by Shi in \cite{Songling1980,Shi1981} and Chen and Wang in \cite{ChenWang1979}, where the authors proved the existence of a quadratic system with 4 limit cycles. For $n=3,$ Li, Liu, and Yang proved in \cite{LiLiuYang2009} that $\mathcal{H}(3)\geq13$. For $n=4$ the authors in \cite{ProTor2019} show  that $\mathcal{H}(4)\geq 28$. Concerning the local version, Bautin showed in \cite{Bautin1954} that $\mathcal M(2)=3$ and,  due to great difficulty this is the only fully resolved case. Following, Giné, Gouveia, and Torregrosa \cite{FerGinTor2021} proved that $\mathcal M(3)\geq12$ and $\mathcal M(4)\geq21$, and Gouveia and Torregrosa \cite{TorreGouv2021} proved $\mathcal M(5)\geq33$.

Regarding Kolmogorov systems, let us define $\mathcal M_{K}(n)$ as the maximal number of limit cycles that can bifurcate from one equilibrium point of the system \eqref{eq:2} with degree $n=m+1$. In this context, we restrict even more the 16th Hilbert Problem since the perturbative parameters must keep the system in the Kolmogorov form. Then the number of limit cycles is expected to be even smaller than $\mathcal M(n).$ In the case $n=1$ we have $\mathcal M_{K}(1)=0$. Bautin proved that $\mathcal M_{K}(2)=0,$ see \cite{Bautin1954}. The case $n=3,$ was studied intensively in the last century \cite{LloydPearson1996}. For a specific family, N. G. Lloyd proved that the maximum order of weak focus is six, and obtained that $\mathcal M_{K}(3)\geq6 $, see \cite{LloydPearson2002}. Regarding $n=4$ the authors shows up a example such that $\mathcal M_{K}(4)\geq7,$ see \cite{HeHaungWang2020}, and in \cite{DuLiuHuang2016} the authors show that $\mathcal M_{K}(4)\geq 8$. Recall that a weak focus is a hyperbolic equilibrium point of type focus with purely imaginary eigenvalues and its order is related to its first non-zero Lyapunov quantity, as we will see later. 

In this work, we will show a new Kolmogorov cubic system with $\mathcal M_{K}(3) \geq 6$. In addition, we improve $\mathcal M_{K}(n)$, for $n=4,5$, we show that $\mathcal M_{K}(4) \geq 13$, and $\mathcal M_{K}(5) \geq 22$. In other words, we give new lower bounds for $\mathcal M_{K}(n)$ for quartic and quintic Kolmogorov systems. Following we show a table with the numbers  $\mathcal M(n)$ and $\mathcal M_{K}(n)$ for $1\leq n\leq 5$.
\begin{table}[h]
\begin{center}
\begin{tabular}{|l|c|c|}
\cline{2-3}
\hline
\multicolumn{1}{|c|}{deg}  & General system & Kolmogorov system   \\
\hline  
$n=1$ & $\mathcal M(1)=0,$  & $\mathcal M_{K}(1)=0,$  \\
\hline 
$n=2$ & $\mathcal M(2)=3,$  & $\mathcal M_{K}(2)=0,$   \\
\hline 
$n=3$ & $\mathcal M(3)\geq 12,$  & $\mathcal M_{K}(3)\geq6,$   \\
\hline 
$n=4$ & $\mathcal M(4)\geq21,$  & $\mathcal M_{K}(4)\geq 13,$  \\ 
\hline
$n=5$ & $\mathcal M(5)\geq33,$  & $\mathcal M_{K}(5)\geq22.$  \\ 
\hline
\end{tabular}
\end{center}
		
\vspace{0.2cm}
		
\caption{Summary of Hilbert numbers for general and Kolmogorov polynomial systems of degree $1\leq n \leq 5.$}
\end{table}

This paper is structured as follows. In section~\ref{se:preliminaries} we show the necessary definitions and algorithms to obtain the coefficients of the return map, the so-called Lyapunov constants. In section~\ref{se:necsuficondition} we dedicate ourselves to the study of center conditions for a family of Kolmogorov systems and we prove the next proposition.

\begin{pro}\label{thmmainn_2} For the cubic Kolmogorov family,
\begin{equation}
\begin{array}{lll}\label{komcub}
\dot{x}&=&x(y-1)(a_2b_1\,x+a_1a_3\,y+a_1b_1-a_1a_3-a_2b_1)/(-b_1^2),\\
\dot{y}&=&y(x-1)(b_1b_2\,x+a_1b_3\,y+a_1b_1-a_1b_3-b_1b_2)/a_1^2,
\end{array}
\end{equation}
the equilibrium point $(1,1)$ is a center if, and only if, one of the next conditions holds:

\begin{itemize}
\item[$(\mathcal{C}_1)$] $a_2=b_2=0;\vspace{0.1cm}$
\item[$(\mathcal{C}_2)$] $a_2=b_3=0;\vspace{0.1cm}$
\item[$(\mathcal{C}_3)$] $a_3=b_2=0;\vspace{0.1cm}$
\item[$(\mathcal{C}_4)$] $a_3=b_3=0;\vspace{0.1cm}$
\item[$(\mathcal{C}_5)$] $a_2-b_2=b_3-a_3=0;\vspace{0.1cm}$
\item[$(\mathcal{C}_6)$] $a_3-b_2=b_1-a_1=b_3-a_2=0;\vspace{0.1cm}$
\item[$(\mathcal{C}_7)$] $a_3+b_2=b_1+a_1=a_2+b_3=0;\vspace{0.1cm}$
\item[$(\mathcal{C}_8)$] $a_1+(a_3^2-b_2^2)/b_2=a_2+b_2=b_1-(a_3^2-b_2^2)/a_3=b_3+a_3=0.\vspace{0.1cm}$
\end{itemize}	
\end{pro}
	
In section~\ref{se:bifurcation}, using Proposition~\ref{thmmainn_2} we show examples of Kolmogorov systems of degree $3, 4$ and $5$ such that $\mathcal M_{K}(3)\geq6$, $\mathcal M_{K}(4)\geq13$, and $\mathcal M_{K}(5)\geq18$. In addition, using the same idea shown in section 3, we will show another example such that $\mathcal M_{K}(4)\geq13$, and $\mathcal M_{K}(5)\geq22$. We summarize these results in the next theorem.

\begin{teo}\label{thmmainn_3} The number of limit cycles bifurcating from a singular monodromic point for a planar Kolmogorov differential system is, at least, $\mathcal M_{K}(3)\geq 6, \mathcal M_{K}(4)\geq 13$, and $\mathcal M_{K}(5)\geq 22$, for a system of degree three, four, and five respectively.
\end{teo}

\section{Preliminaries}\label{se:preliminaries}
	
In this section, we recall how to obtain the coefficients of the return map near an equilibrium point of monodromic or center-type, a weak focus, or a center, the so-called Lyapunov Constants. Let us consider a polynomial differential system of degree $n$ having an equilibrium point at the origin such that the eigenvalues of their Jacobian matrix are purely imaginary. For simplicity, we will consider only the cases when the linear part of each system is written in its normal form as follows
\begin{equation}\label{eq:3_1}
\begin{array}{lll}
\dot{x}&=&-y + \displaystyle\sum_{k=2}^n g_k(x,y),\\
\dot{y}&=&\phantom{-}x+\displaystyle\sum_{k=2}^n f_k(x,y),
\end{array}
\end{equation}
where  $g_k,$ $f_k$ denote homogeneous polynomials of degree $k$ for $n \geq k\geq 2$ such that
\begin{equation*}
	g_k(x,y)=\sum_{\ell=0 }^k a_{k-\ell,\ell}x^{kp-\ell}y^\ell \; \mbox{and} \; 	f_k(x,y)=\sum_{\ell=0 }^k b_{k-\ell,\ell}x^{kp-\ell}y^\ell.
\end{equation*}

\subsection{The center conditions}\label{se:centerconditions}

Given a planar polynomial differential system with a critical point at the origin whose linearization provides a center, under what conditions can we conclude that the origin is a center for the nonlinear system? If the determinant of the linear part at the origin is non-zero then the origin is non-degenerate. Therefore, as long as the linear parts do not give a center, i.e, the eigenvalues of the linear part are not pure imaginary then the Hartman–Grossman theorem tells us that in a sufficiently small neighborhood of the origin, the system is topologically equivalent to its linear part, i.e., we can ignore the terms of higher order. In our situation, the center-focus problem asks for the criteria which determine whether the origin whose linear parts give a center, really is a center. If the origin is either a center or focus, we shall use the more general term monodromic to cover both cases. A focus whose linearization gives a center is called a weak focus. 	

The approach to characterize when system~\eqref{eq:3_1} has a center at the origin is based on the well-known Poincar\'e--Lyapunov Theorem, see \cite{Liapounoff1965,Poincare1886}. A non-constant analytical function defined in a neighborhood $\Omega$ of the origin, $\Phi:\Omega\subset\mathbb{R}^2\rightarrow \mathbb{R}^2,$ is a first integral of system \eqref{eq:3_1} if it is constant along any solution $\gamma$ or, equivalently,
\begin{equation}\label{eq:3}
\left.\dot{x}\,\frac{\partial \Phi}{\partial x}+\dot{y}\,\frac{\partial \Phi}{\partial y}\right|_{\gamma}\equiv 0.
\end{equation}
It is possible characterize centers in terms of their first integrals using the well-known Poincaré-Lyapunov Center Theorem, see \cite{IliYak2008,RomaSha2009}. The system~\eqref{eq:3_1} has a center at the origin if and only if it admits a local analytic first integral of the form
		\begin{equation}\label{eq:2_1}
			\Phi(x,y)=x^2+y^2 + \sum_{p=3}^\infty F_p(x,y), \; \mbox{where} \; F_p(x,y)=\sum_{\ell=0 }^p q_{p-\ell,\ell}x^{p-\ell}y^\ell.
		\end{equation}
In addition, the existence of a formal first integral $\Phi$ of the above form implies the existence of a local analytic first integral.

The necessary conditions for the existence of a first integral \eqref{eq:2_1} for system \eqref{eq:3_1} are obtained by looking for a formal series \eqref{eq:2_1} satisfying \eqref{eq:3}. Although \eqref{eq:3} is not always satisfied, it is always possible to choose coefficients of the formal power series \eqref{eq:2_1} which satisfies the following equation:
\begin{equation}\label{eqliapu}
	\dot{x}\,\frac{\partial \Phi}{\partial x}+\dot{y}\,\frac{\partial \Phi}{\partial y}=\sum_{j=1}^\infty L_j (x^2+y^2)^{j+1},
\end{equation}
see \cite{BlowLlo1984,Lia1966,NemSte1960,Shi1981}. Any non-zero $L_j$ is an obstruction for the origin to be a center. Furthermore, the quantities $L_j$ are rational functions whose numerators are polynomials depending on the parameters of system~\eqref{eq:3_1}. 

In general in \eqref{eqliapu}, the coefficients of the homogeneous terms of an odd degree have a unique solution, since the determinant of its associated matrix never vanishes. In contrast, those of even degree does not have a single solution, after all the determinant always vanishes. The quantities $L_j$ ensure that linear systems corresponding to even degrees are possible and indeterminate. Therefore, the different choice of solutions means that the next quantity $L_{j+1}$ will be changed. Considering the ring of homogeneous polynomials in coefficients of system~\eqref{eq:3_1}, the authors in \cite{Shi1984} proved that the quantity $L_j$, of order $j\geq1$ is uniquely determined modulo the ideal generated by $L_1, \dots, L_{j-1} $ for $L_0=constant$.

Starting the computational procedure for finding the first $N$ conditions for integrability, we write down \eqref{eq:2_1} up to order $2N+2$
\begin{equation}\label{eq:4}
\widetilde{\Phi}(x,y)=\dfrac{x^2+y^2}{2}+\sum_{p= 3}^{2N+2}\sum_{\ell=0 }^p q_{p-\ell,\ell}x^{p-\ell}y^\ell.
\end{equation}
Then, for each $i=3,\ldots,2N+2$ we equate to zero the coefficients of terms of degree $i$ in the expression 
\begin{equation*}
\dot{x}\frac{\partial \widetilde\Phi}{\partial x}+\dot{y}\frac{\partial \widetilde\Phi}{\partial y}=\left(-y+\displaystyle\sum_{k=2}^n g_k(x,y)\right)\frac{\partial \widetilde\Phi}{\partial x}+\left(x+\displaystyle\sum_{k=2}^n f_k(x,y)\right)\frac{\partial \widetilde\Phi}{\partial y}.
\end{equation*}

Starting at $i=3,$ we should solve in a recurrent way each linear system of $i+1$ equations with $i+1$ variables, $q_{p-\ell,\ell}$ such that $\ell=0,\dots,p$. All linear systems corresponding to odd degrees, $i=2j+1,$ have a unique solution in terms of previous values of $q_{p-\ell,\ell}$. As the determinant of the linear system that corresponds to an even degree, $i=2j+2,$ vanishes, we need to add an extra condition so that the linear system has a unique solution. In fact, at this step, we have one equation more than the number of variables. We add suitable equations, for the terms $(x^{2}+y^{2})^{j+2}$ for example, so that the derivative over the associated vector field writes as
\begin{equation}\label{eq:6}
\dot{x}\frac{\partial\Phi}{\partial x}+\dot{y}\frac{\partial \Phi}{\partial y}=\sum_{j=1}^\infty L_{j}(x^{2}+y^{2})^{j+2}.
\end{equation}

When at least one $L_{j}$ is different from zero, follow that $\Phi$ is a Lyapunov function in a neighborhood of the origin. Then, system~\eqref{eq:3_1} has no local analytic first integral and we say that the equilibrium point is a \emph{weak focus} of order $k$ if the first nonzero coefficient in \eqref{eq:6} is $L_k$. The coefficient $L_j$ in \eqref{eq:6} is called the $j$-th Lyapunov constants or focus value. Furthermore, the first nonzero quantity $L_j$ and its order are invariants \cite{Shi1984}. 

Denoting by $\lambda \in \mathbb{R}^M$ the set of parameters of the system~\eqref{eq:3_1} the $L_j$'s will be polynomials in $\lambda$. In addition, the set $\mathcal{B}^{\mathbb{R}}=\langle L_1,L_2,\ldots\rangle$  will be an ideal in the polynomial ring $\mathbb{\mathbb{R}}[\lambda]$. The real importance of ideal $\mathcal{B}^{\mathbb{R}}$ follows that if all $L_j$ vanishes, then all Lyapunov quantities will vanish regardless of the way to calculate them. 
In this aspect, we introduce the next definition which recalls the notion of \emph{Bautin ideal} and the \emph{center variety}.

\begin{defini}
The ideal defined by the real focus quantities, $\mathcal{B}^{\mathbb{R}}=\langle L_1,L_2,\ldots\rangle\subset\mathbb{R}[\lambda]$ is called the real Bautin ideal. The affine variety $ \mathbf{V}^{\mathbb{R}}=\mathbf{V}(\mathcal{B}^\mathbb{R})$ is called the real center variety of system~\eqref{eq:1}, such that
\begin{equation*} 
\mathbf{V}(\mathcal{B}^{\mathbb{R}})=\{\lambda \in \mathbb{R}^M:r(\lambda)=0, \; \forall \; r \in \mathcal{B}^{\mathbb{R}}\}.
	\end{equation*}
\end{defini}

When we can explicitly determine the real center variety we have the center-focus problem solved for the system \eqref{eq:1}. However, in most cases, this is not a simple problem. On the other hand, the Hilbert Basis Theorem \cite{AdaLou1994,CoxLitOsh1997,MacBir1967} assures that $V(B)$ is finitely generated. Then there exists a positive integer $j$ such that $\mathcal{B}^{\mathbb{R}}=\mathcal{B}^{\mathbb{R}}_j=\langle L_1,\ldots,L_j\rangle.$ In other words, we know that for $j$ big enough, the above algorithm provides a necessary set of conditions, $\{L_j=0: j=1,\ldots, N\}$, for system~\eqref{eq:3_1} be a center. The main difficulty follow that there is no technique to obtain $j$ a priori. We can also say that the polynomials $L_j$ represent obstacles to the existence of a first integral. In particular, system~\eqref{eq:3_1} admits a first integral of the form~\eqref{eq:2_1} if and only if $L_j=0,$ for all $j\geq1.$ Thus, the simultaneous vanishing of all focus quantities provides conditions that characterize when a system of the form~\eqref{eq:3_1} has a center at the origin.

Notice that the inclusion $\mathbf{V}^{\mathbb{R}}=\mathbf{V}(\mathcal{B}^\mathbb{R})\subset\mathbf{V}(\mathcal{B}^\mathbb{R}_j)$ holds for any $j\geq1.$ The opposite inclusion, for a fixed $j$, is verified by finding the irreducible decomposition of $\mathbf{V}(\mathcal{B}^\mathbb{R}_j),$ see \cite{RomaSha2009}, such that any point of each component of the decomposition corresponds to a system having a center at the origin. 

At this point, the importance of Lyapunov constants in the center and cyclicity problems has been made clear.

\begin{remark}
	There are many other methods of calculating equivalent sets of quantities $L_j$. For example, the same analysis described previously but considering the vector field in complex coordinates, via the change of variables $z = x + i\,y$, see \cite{GarMazSha2018,LiuHua2005,RomaSha2009}. Another interesting technique is the Andronov method \cite{AndLeoGorMai1973}, which is actually the classical technique to obtain equivalent expressions of this quantities.
\end{remark}

\subsection{Degenerate Hopf bifurcation}

We are interested only in the limit cycles bifurcating from equilibrium points of the center or focus type, where the return map is well-defined and analytic. The classical Hopf bifurcation, see \cite{MarMcC1976}, occurs when a complex conjugate pair of eigenvalues of the linearised system at an equilibrium point becomes purely imaginary. Then the birth of a limit cycle from a weak focus of first-order arises varying the trace of the matrix of the linearized system from zero to a number small enough. In this way, we will denote by $L_0$ the trace of the perturbed system, clearly $L_0=0$ for the unperturbed one. The degenerate Hopf bifurcation is the natural generalization of this bifurcation phenomenon when $k$ small limit cycles appear from a weak focus of order $k$. In general, the complete unfolding of $k$ limit cycles near a weak focus of order $k$ is only guaranteed when we analyze for analytic perturbations, see for example \cite{Rou1998}. When the perturbation is restricted to be a polynomial of some fixed degree this property is not automatic. Because of this, the problem of finding the cyclicity of a center, for systems like \eqref{eq:3_1}, is so difficult. A way to avoid these difficulties is to study lower bounds for cyclicity. In \cite{Chr2006}, Christopher provides the necessary conditions to obtain lower bounds for the cyclicity of a center. Christopher details how we can use the first-order Taylor approximation of the Lyapunov constants to obtain lower bounds on the number of limit cycles near center-type equilibrium points. There are similar previous results due to Chicone and Jacobs (\cite{ChiJac1989,ChiJac1991}), and due to Han \cite{Han1999}.

\begin{pro}[\cite{Chr2006}]\label{thm:chrorder1}
Suppose that $s$ is a point on the center variety and that the first $k$ Lyapunov constants, $L_{1},\ldots,L_{k},$ have independent linear parts (with respect to the expansion of $L_{i}$ about $s$), then $s$ lies on a component of the center variety of codimension at least $k$ and there are bifurcations which produce $k$ limit cycles locally from the center-type equilibrium point corresponding to the parameter value $s$. If, furthermore, we know that $s$ lies on a component of the center variety of codimension $k$, then $s$ is a smooth point of the variety, and the cyclicity of the center for the parameter value $s$ is exactly $k$. In the latter case, $k$ is also the cyclicity of a generic point on this component of the center variety.
\end{pro}

When the linear parts of the next Lyapunov constants are a linear combination of the firsts $k$-th we can use the higher developments to obtain more limit cycles. The next result, also due to Christopher in \cite{Chr2006}, is an extension of Proposition~\ref{thm:chrorder1} that shows as sometimes we can obtain more limit cycles using high-order Taylor developments of the Lyapunov constants. 

\begin{pro}[\cite{Chr2006}]
Suppose that we are in a point $s$ where Proposition~\ref{thm:chrorder1} applies. After a change of variables if necessary, we can assume that $L_0=L_1=\cdots=L_k=0$ and the next Lyapunov constants $L_{i}=h_i(\lambda)+O_{m+1}(\lambda),$ for $i=k+1,\ldots,k+l,$ where $h_i$ are homogeneous polynomials of degree $m\ge2$ and $\lambda=(\lambda_{k+1},\ldots,\lambda_{k+l})$. If there exists a line $\boldsymbol \ell,$ in the parameter space, such that $h_{i}(\boldsymbol \ell) = 0,$ $i=k+1,\ldots,k+l-1,$ the hypersurfaces $h_{i} = 0$ intersect transversally along $\ell$ for $i=k+1,\dots,k+l-1,$ and $h_{k+l}(\boldsymbol \ell)\ne 0$, then there are perturbations of the center which can produce $k+l$ limit cycles.
\end{pro}

The above result is not written exactly as in the original Christopher paper because we have adapted to include also the conclusion of Proposition~\ref{thm:chrorder1}. And it is important to mention that we have the number of the parameters $k+l=M$ or we consider the parameters with index $k+l+1,\dots, M$ being nulls. 

We observe that Proposition~\ref{thm:chrorder1} is devoted to obtaining limit cycles using the linear terms of the Lyapunov Constants but for fixed systems. Then natural questions arise is if the same method is can be used for family systems. The answer is positive and the next theorem whose proof can be found in \cite{FerGinTor2021} extend Proposition~\ref{thm:chrorder1}. Let $(\dot{x},\dot y)=(P_c(x,y,\mu),Q_c(x,y,\mu))$ be a family of polynomial centers of degree $n$ depending on a parameter $\mu\in\mathbb{R}^\ell,$ having a center equilibrium point at the origin. Then we consider the perturbed polynomial system
\begin{equation}\label{eq:centrob}
\begin{array}{lll}
\dot{x} & =&P_c(x,y,\mu)+\alpha x+P(x,y,\lambda),\\
\dot{y} & =&Q_c(x,y,\mu)+\alpha y+Q(x,y,\lambda),\\
\end{array}
\end{equation}	
with $P,Q$ polynomials of degree $n$ having no constant nor linear terms. More concretely,
\[P(x,y,\lambda)=\sum\limits_{k+l=2}^na_{k,l}x^ky^l,\quad Q(x,y,\lambda)=\sum\limits_{k+l=2}^nb_{k,l}x^ky^l,\]
with $\lambda=(a_{20},a_{11},a_{02},\ldots,b_{20},b_{11},b_{02})\in\mathbb{R}^{n^2+3n-4}= \mathbb{R}^{M}.$ And, the trace parameter $2 \, \alpha$ sometimes is also denoted by $L_0$, as mentioned before.

\begin{pro}\label{theop}
When $\alpha=0,$ we denote by $L_{j}^{(1)}(\lambda,\mu)$ the first-order development, with respect to $\lambda\in \mathbb{R}^{M},$ of the $j$-Lyapunov constant of system~\eqref{eq:centrob}. We assume that, after a change of variables in the parameter space if necessary, we can write

\begin{equation*}\label{eqL}
	L_j=\begin{cases}
				&\lambda_j + O_2(\lambda), \text{ for } j=1,\ldots,k-1,\\
				&\sum\limits_{l=1}^{k-1} g_{j,l}(\mu) \lambda_l+f_{j-k}(\mu)\lambda_{k}+ O_2(\lambda), \text{ for } j=k,\ldots,k+\ell.
			\end{cases}
\end{equation*}
Where with $O_2(\lambda)$ we denote all the monomials of degree higher or equal than $2$ in $\lambda$ with coefficients analytic functions in $\mu$. If there exists a point $\mu^*$ such that $f_0(\mu^*)=\cdots=f_{\ell-1}(\mu^*)=0,$ $f_{\ell}(\mu^*)\ne0,$ and the Jacobian matrix of $(f_{0},\ldots,f_{\ell-1})$ with respect to $\mu$ has rank $\ell$ at $\mu^*,$ then system~\eqref{eq:centrob} has $k+\ell$ hyperbolic limit cycles of small amplitude bifurcating from the origin. 
\end{pro}
	
\begin{remark}
We remark on the importance of the number of components in parameters $\lambda$ and $\mu$ in the theorem above. If there are more parameters than the relevant $k$ in $\lambda$, in $O_2(\lambda)$ term can appear monomials of degree $2$ that can affect the monomials of degree $1$ and the result could be not valid. The above result shows that on some special points on such components, the cyclicity can increase. We use this mechanism to obtain another cubic Kolmogorov system with six limit cycles and to improve the lower bound of limit cycles for the quartic Kolmogorov system. 
\end{remark}
	
\section{Centers classification}\label{se:necsuficondition}	

We devote this section to proving the center classification statements and the necessary conditions of Proposition~\ref{thmmainn_2}. As the proof is quite long we have written it separately in two propositions, Proposition~\ref{cent1} and  Proposition~\ref{cent2}.
	
As in general we do not know how many focal coefficients are sufficient to determine the Bautin variety, the solving of the Center-Focus Problem for a specific system is divided into steps: in the first we calculate a finite number of focal values. After that, the second step is devoted to determining what are the conditions under the system parameters that cancel these focal coefficients. The third step consists of verifying if the conditions found are sufficient for us to guarantee that the origin is a center. In general, this last step involves either the application of Darboux's Theory of Integrability, see \cite{DumLliArt2006}, or the investigation of symmetries or reversibility of the system \cite{RomaSha2009}.  

\begin{pro}\label{cent1} If the equilibrium point $(1,1)$ of the Kolmogorov cubic system \eqref{komcub} is a center, then the parameters $(a_1,a_2,a_3,b_1,b_2,b_3)$ satisfy one of the conditions given in the statement of Proposition~\ref{thmmainn_2}.      
\end{pro}
\begin{proof} 
Taking the equilibrium point $(1,1)$ of the system~\eqref{komcub}, follow that this singular point is a nondegenerate equilibrium point of center-focus type, that is, the trace and the determinant of the Jacobian matrix are zero and a positive number, respectively. As usual, we need to compute a few Lyapunov Constants. As we have six parameters, we must compute, at least, six Lyapunov Constants. Then, we have a system with six equations and six parameters
\begin{equation*}\label{S}
\mathcal{S}=\{L_1=L_2=L_3=L_4=L_5=L_6=0\}.
\end{equation*}
By solving the obtained system, we can check that all the solutions correspond to the ones described in the statement. Now we do an affine change of coordinates given by
\begin{equation*}\label{change}
\left(u,v\right)\rightarrow \left((x-1)/a_1,(y-1)/b_1\right).
\end{equation*} 
Then we translate the point $(1,1),$ to the origin and we have
\begin{equation*}
\begin{array}{lll}\label{komcub_chang}
\dot{u}&=&-v(a_1\,u+1)(a_2\,u+a_3\,v+1),\\
\dot{v}&=&\phantom{-}u(b_1\,v+1)(b_2\,u+b_3\,v+1).
\end{array}
\end{equation*}
Following the approach described in Section~\ref{se:centerconditions} for the computation of the center conditions $L_i$, $i=1,\ldots,6$, we can obtain the six first Lyapunov constants that are polynomial in the parameters $(a_1,a_2,a_3,b_1,b_2,b_3).$ Due to the size of expressions, we will show only the first three Lyapunov Constants. 
\begin{equation*}
\begin{aligned}
L_1=&\phantom{-}\frac{1}{3}(a_3a_2-b_2b_3),\\
L_2=&-\frac{b_2b_3}{15}(2a_1a_2 - 2a_1b_2 - a_2^2 - 4a_2b_2 - 5a_3^2 + 2a_3b_1 + 4a_3b_3 - 2b_1b_3 + 5b_2^2 + b_3^2),\\
L_3=&-\frac{1}{210}(b_3b_2(2a_1a_2^3 - 14a_1a_2^2b_2 + 2a_1a_2b_2^2 - 2a_1a_2b_3^2 + 10a_1a_3^2b_2 - 22a_1a_3b_2b_3 \\& - 105a_1b_2^3 + 14a_1b_2b_3^2 - a_2^4 - 4a_2^3b_2 + 30a_2^2b_2^2 + a_2^2b_3^2 - 20a_2b_2^3 - 2a_2b_2b_3^2 + 30a_3^4 \\&- 30a_3^3b_3 - 25a_3^2b_2^2 - 6a_3^2b_3^2 + 50a_3b_2^2b_3 + 6a_3b_3^3 - 5b_2^4 - 24b_2^2b_3^2)).
\end{aligned}
\end{equation*}

Now we need to solve the algebraic system of equations \eqref{S}, however, despite this system having only six variables and equations the usual mechanisms for solving it fail. Then, to determine the irreducible components of the variety 
\begin{equation*}\label{V}
\textbf{V}=\textbf{V}(L_1,L_2,L_3,L_4,L_5,L_6),
\end{equation*}
consists of using resultants, see \cite{Sturm2002}. We will use the crossover of the resultants of the Lyapunov Constants constants $L_i$ with $i=1,\ldots,6,$ regarding of variables. Computing the five resultants concerning the parameter $a_2.$ we have
\begin{equation*}
\begin{aligned}
\Res(L_1,L_2,a_2)=&\mathcal{R}_1\cdot \mathcal{R}_2\cdot \mathcal{R}_{12} , \quad \Res(L_1,L_3,a_2)=\mathcal{R}_1\cdot \mathcal{R}_2\cdot\mathcal{R}_{13},\\
\Res(L_1,L_4,a_2)=&\mathcal{R}_1\cdot \mathcal{R}_2\cdot\mathcal{R}_{14},\quad \Res(L_1,L_5,a_2)=\mathcal{R}_1\cdot \mathcal{R}_2\cdot \mathcal{R}_{15},\\
\Res(L_1,L_6,a_2)= & \mathcal{R}_1\cdot \mathcal{R}_2\cdot\mathcal{R}_{16},
\end{aligned}
\end{equation*}  
where $\mathcal{R}_1=b_2,$ $\mathcal{R}_2=b_3,$ and $\mathcal{R}_{1j}$ for $j=2,\ldots,6,$ are polynomials of degrees $4,8,11,15,13$ in $(a_1,a_3,b_1,b_2,b_3),$ respectively. In a second step, we compute resultant of $\mathcal{R}_{1j}$ for $j=2,\ldots,6,$ with respect to the parameter $a_1.$ Then we have 
\begin{equation*}
\begin{aligned}
\Res(\mathcal{R}_{12},\mathcal{R}_{13},a_1)=& \mathcal{R}_{123} , \quad \Res(\mathcal{R}_{12},\mathcal{R}_{14},a_1)= \mathcal{R}_{124},\\
\Res(\mathcal{R}_{12},\mathcal{R}_{15},a_1)=& \mathcal{R}_{125},\quad
\Res(\mathcal{R}_{12},\mathcal{R}_{16},a_1)= \mathcal{R}_{126},\\
\end{aligned}
\end{equation*}   
where $\mathcal{R}_{12j}$ for $j=3,\ldots,6,$ are polynomials of degrees $11,17,18,13,$ in $(a_3,b_1,b_2,b_3),$ respectively.  Following the same scheme we calculate the resultant using $\mathcal{R}_{12j}$ with respect to the parameters $b_1$, and we obtain 
\begin{equation*}
\begin{aligned}
\Res(\mathcal{R}_{123},\mathcal{R}_{124},b_1)=&\mathcal{R}_1^4\cdot\mathcal{R}_3^8\cdot\mathcal{R}_4^3\cdot\mathcal{R}_5^3\cdot\mathcal{R}_6\cdot \mathcal{R}_{1234},\\
\Res(\mathcal{R}_{123},\mathcal{R}_{125},b_1)=&\mathcal{R}_1^2\cdot\mathcal{R}_3^4\cdot\mathcal{R}_4^2\cdot\mathcal{R}_5^2\cdot\mathcal{R}_6\cdot \mathcal{R}_{1235},\\
\Res(\mathcal{R}_{123},\mathcal{R}_{126},b_1)=&\mathcal{R}_3\cdot\mathcal{R}_4\cdot\mathcal{R}_5\cdot\mathcal{R}_6\cdot \mathcal{R}_{1236},
\end{aligned}
\end{equation*}  
where $\mathcal{R}_3=a_3-b_3,$ $\mathcal{R}_4=a_3-b_2,$ $\mathcal{R}_5=a_3+b_2,$ $\mathcal{R}_6=a_3+b_3,$ and $\mathcal{R}_{123j}$ with $j=4,\ldots,6,$ are polynomials of degrees $18, 17, 9$ in $(a_3,b_2,b_3),$ respectively. Finally, we calculate the resultant using $\mathcal{R}_{123j}$ with respect to the parameters $b_2,$ and we obtain
\begin{equation*}
\begin{aligned}
\Res(\mathcal{R}_{1234},\mathcal{R}_{1235},b_2)=&\mathcal{R}_7^{64}\cdot \mathcal{R}_{12345}, \quad \Res(\mathcal{R}_{1234},\mathcal{R}_{1236},b_2)=&\mathcal{R}_7^{64}\cdot \mathcal{R}_{12346},
\end{aligned}
\end{equation*}   
where $\mathcal{R}_7=a_3$ and $\mathcal{R}_{1234j}$ for $j=5,6$ are polynomials of degrees $88, 84,$ in $(a_3,b_3),$ respectively. In addition, solving the algebraic system $\{\mathcal{R}_{12345}=\mathcal{R}_{12346}=0\},$ calling for $\mathcal{R}_8=a_3+b_3/5,$ and $\mathcal{R}_9=a_3-b_3/3,$ we obtain the solutions $\mathcal{R}_i=0,$ with $i=7,8,9.$ Then, to find the irreducible components of the variety \eqref{V}, we can branch the algebraic system \eqref{S}, into the nine following algebraic systems:
\begin{equation*}
	\mathcal{S}_i=\{L_1=L_2=L_3=L_4=L_5=L_6=\mathcal{R}_i=0\}, \quad i=1,\ldots,9.
\end{equation*}
Thus, from the nine algebraic systems above we obtain twenty-one solutions, and by filtering these solutions we obtain all the centers given in the statement.
\end{proof}

\begin{pro}\label{cent2} For each family $\mathcal{C}_i,$ $i=1,..,8$, listed in Proposition~\ref{thmmainn_2}, the corresponding Kolmogorov system \eqref{komcub} has a center at the equilibrium point $(1,1).$       
\end{pro}

\begin{proof} 
We start to do an affine change of coordinates \eqref{change} that allows us to translate the point $(1,1),$ to the origin and obtain the system in the form \eqref{komcub_chang}. We observe that in this case, the origin is a nondegenerate center-focus point. The proof follows straightforwardly doing a case-by-case study. Using the Darboux theory to compute the integration factor and reversibility properties we show the existence of a center at the origin. 
We observe that the cases $\mathcal{C}_{6}$) and $\mathcal{C}_{7}$) we have the reversibility in relation of straight line $y=x.$ Following we show the corresponding integration factor.
\begin{equation*}
	\begin{aligned}
		\mathcal{C}_{1}): 
		V & =\frac{a_1\,b_1^3\,b_3^3}{(b_1\,y+1)(a_1\,x+1)(b_3\,y+1)}, \; \; \; \; \; \; \;
		\mathcal{C}_{2}):  V = \frac{a_1\,b_1}{(a_1\,x+1)(b_1\,y+1)},\\
		\mathcal{C}_{3}):   V & =\left(\frac{b_1\,y+1)\,(a_1\,x+1)}{a_1\,b_1}\right)^{\alpha}\left(\frac{(b_3\,y+1)\,(a_1\,x+1)}{a_1\,b_3}\right)^{\beta}\left(\frac{(a_1\,x+1)}{a_1}\right)^{\gamma}, \\
			\mathcal{C}_{4}):  V & =\frac{a_1\,b_1\,a_2}{(b_1\,y+1)\,(a_1\,x+1)\,(a_2\,x+1)}, \\
			\mathcal{C}_{5}):  V &  = \frac{b_1\,b_2\,a_1}{2(b_2\,x+b_3\,y+1)\,(b_1\,y+1)\,(a_1\,x+1)}, \\
			\mathcal{C}_{8}): V & =\frac{\delta\,(a_3\,y+b_2\,x+1)}{(a_3^2\,y-b_2^2\,y+a_3)\,(a_3^2\,x-b_2^2\,x-b_2)\,(a_3^4\,y^2-a_3^2\,b_2^2\,x^2-a_3^2\,b_2^2\,y^2+b_2^4\,x^2-a_3^2-b_2^2)},
		\end{aligned}
\end{equation*}
\smallskip	
\noindent
where 
\begin{equation*}
	\begin{aligned}
		\alpha & = -\frac{a_1\,a_2-a_2^2+b_1^2-b_1\,b_3}{(b_1-b_3)\,b_1}, \; \; \; \; \; \; \; \; \; \; \; \; \; \; \; \; \; \; \; \; \beta  = \frac{a_1\,a_2-a_2^2-b_1\,b_3+b_3^2}{b_3\,(b_1-b_3)}, \\
		\gamma & = -\frac{a_1^2\,a_2-a_1\,a_2^2-a_1\,b_1\,b_3+a_2\,b_1\,b_3}{a_1\,b_1\,b_3}, \;	\; \; \; \;  \delta =  a_3\,b_2\,(a_3^2+b_2^2).
	\end{aligned}
\end{equation*}
\end{proof}
	
\section{Bifurcation of limit cycles in Kolmogorov systems}\label{se:bifurcation}	
	
This section is devoted to proofing the Theorem~\ref{thmmainn_3}. The prove of theorem follow from the next Propositions~[\ref{pro6cl},\ref{pro13cl1},\ref{pro13lc2} and \ref{pro22cl}].  For the following result, we use one among the families of cubic Kolmogorov centers from Proposition~\ref{thmmainn_2}. Doing a brief analysis, we selected the most promising family to obtain limit cycles after adequate perturbations.

\begin{pro}\label{pro6cl} There exist parameters perturbations such that the center family $\mathcal{C}_8,$ given by 
\begin{equation}
\begin{array}{lll}\label{komcub_c1}
\dot{x}&=&x(y-1)(b_2^2\,x+a_3^2\,y-2b_2^2)a_3/\left[(a_3-b_2)(a_3+b_2)b_2\right],\\
\dot{y}&=&y(x-1)(b_2^2\,x+a_3^2\,y-2a_3^2)b_2/\left[(a_3-b_2)(a_3+b_2)a_3\right],
\end{array}
\end{equation}	
can produce, at least, six small limit cycles bifurcation from $(1,1).$     
\end{pro}
\begin{proof}[Proof]
Firstly, we will consider $b_ 2=2$, and we will make a change of coordinates such that the origin is a center. Then, the unperturbed system is given by
\begin{equation*}
	\begin{array}{lll}\label{komcub_c2}
		\dot{x}&=&-(a_3\,y - 2\,x + 1)[1-(a_3^2 - 4)\,x/2]\,y,\\
		\dot{y}&=&\phantom{-}(2\,x-a_3\,y + 1)[1+(a_3^2 - 4)\,y/a_3]\,x,
	\end{array}
\end{equation*}
Now, considering the perturbed system we have 
\begin{equation*}
	\begin{array}{lll}
		\dot{x}=&-(a_3\,y - 2\,x + 1)[1-(a_3^2 - 4)\,x/2]\,y+[1-(a_3^2 - 4)\,x/2](a_{02}\,y^2 + a_{11}\,x\,y + a_{20}\,x^2),\\
		\dot{y}=&\phantom{-}(2\,x-a_3\,y + 1)[1+(a_3^2 - 4)\,y/a_3]\,x+[1+(a_3^2 - 4)\,y/a_3](b_{02}\,y^2 + b_{11}\,x\,y + b_{20}\,x^2).
	\end{array}
\end{equation*}

Computing the linear terms of the first Lyapunov Constants in relation with perturbative parameters $\Lambda=\left\{a_{20}, a_{11}, a_{02}, b_{20}, b_{11}, b_{02} \right\}$ we obtain that the rank of Jacobian Matrix in relation to $\Lambda$ is five. Due to size, we will show only the first three Lyapunov Constants.
\begin{equation*}
\begin{aligned}
L_1= & \phantom{-}\frac{2}{3a_{3}}(-a_{20}a_{3}^3 - a_{11}a_{3}^2 + a_{3}^2b_{02} + a_{3}^2b_{20} + 2a_{02}a_{3} + 2a_{20}a_{3} - 2a_{3}b_{11} - 8b_{02}),\\
L_2= & -\frac{2}{45a_{3}}(-6a_{20}a_{3}^7 - 6a_{11}a_{3}^6 + 6a_{3}^6b_{02} + 6a_{3}^6b_{20} + 12a_{02}a_{3}^5 - 35a_{20}a_{3}^5 - 12a_{3}^5b_{11}  50a_{11}a_{3}^4 \\ & + 8a_{3}^4b_{02} + 38a_{3}^4b_{20} + 16a_{02}a_{3}^3 - 160a_{20}a_{3}^3 - 88a_{3}^3b_{11} - 176a_{11}a_{3}^2 - 320a_{3}^2b_{02} + 32a_{3}^2b_{20} \\ & + 304a_{02}a_{3} + 448a_{20}a_{3} - 400a_{3}b_{11} - 1504b_{02}), \\
L_3 = & \frac{1}{15120a_3^3}(-2253a_{20}a_{3}^{13} - 2253a_{{{11}}}a_{3}^{12} + 2253a_{3}^{12}b_{02} + 2253a_{3}^{12}b_{20} + 4506a_{02}a_{3}^{{11}} \\ & - 41066a_{20}a_{3}^{{11}} - 4506a_{3}^{11}b_{11} 
- 46772a_{11}a_{3}^{10} + 31148a_{3}^{10}b_{02} + 41972a_{3}^{10}b_{20}  \\ & + 59944a_{02}a_{3}^9 - 148952a_{20}a_{3}^9 - 88744a_{3}^9b_{11} - 219648a_{11}a_{3}^8 
- 165664a_{3}^8b_{02}  \\ &+ 122496a_{3}^8b_{20} + 147840a_{02}a_{3}^7 - 423168a_{20}a_{3}^7 - 409728a_{3}^7b_{11} - 645504a_{11}a_{3}^6  
\\ &- 1372800a_{3}^6b_{02} + 153984a_{3}^6b_{20} + 599808a_{02}a_{3}^5 - 545792a_{20}a_{3}^5 - {13}49376a_{3}^5b_{11}  \\ & - 1468160a_{11}a_{3}^4- 4764928a_{3}^4b_{02} + 380672a_{3}^4b_{20}+ 2684416a_{02}a_{3}^3 + 4{12}5184a_{20}a_{3}^3  \\ & - 3483{13}6a_{3}^3b_{11} + 3072a_{11}a_{3}^2 - 13026304a_{3}^2b_{02} - 3072a_{3}^2b_{20} - 6144a_{02}a_{3} - 6144a_{20}a_{3} \\ & + 6144a_{3}b_{11} + 24576b_{02}).
\end{aligned}		
\end{equation*}

Taking the first four Lyapunov Constants with the parameters $\left\{a_{02}, a_{11}, a_{20}, b_{02}\right\}$ we have the Jacobian Matrix has rank four. Then make a change of coordinates, we can write 
\begin{equation*}
\begin{aligned}
L_1= & u_{1}+O_{2}(u_1,u_2,u_3,u_4,a_3,b_{20}), \; \; \; \; \; \; \; \; \; \; \; \; \; \; \, L_2=  u_{2}+O_{2}(u_1,u_2,u_3,u_4,a_3,b_{20}), \\
L_3= & u_{3}+O_{2}(u_1,u_2,u_3,u_4,a_3,b_{20}),  \; \; \; \; \; \; \; \; \; \; \; \; \; \; \, L_4=  u_{4}+O_{2}(u_1,u_2,u_3,u_4,a_3,b_{20}), \\
L_5= &f_{1}(a_3,b_{20})+O_{2}(u_1,u_2,u_3,u_4,a_3,b_{20}), \; \; \;  L_6=  f_{2}(a_3,b_{20})+O_{2}(u_1,u_2,u_3,u_4,a_3,b_{20}), \\
\end{aligned}		
\end{equation*}
 where $f_{1}(a_3,b_{20})$ and $f_{2}(a_3,b_{20})$ are rational polynomials. The numerator part of $f_{1}$ has degree twenty-one in $a_3$ and has twenty monomials. The denominator part of $f_{1}$ has degree twelve in $a_3$ and has seven monomials. Regarding $f_{2}$, the numerator part has degree twenty-five and twelve monomials, and the denominator part has degree twenty-four and seven monomials. Moreover, there is a point $\widetilde{a_3}$ such that $f_{1}(\widetilde{a_3},b_{20})=0$, $f_{2}(\widetilde{a_3},b_{20})\neq 0$, and $f^\prime_{1}(\widetilde{a_3},b_{20}) \neq 0$. Therefore, by Theorem~\ref{theop}, we can obtain six small limit cycles bifurcating from the origin.
\end{proof}


\begin{remark}
We remark that the system shown in Proposition~\ref{pro6cl} is the second example in the literature in which 6 small limit cycles bifurcating from a center. As mentioned, the first example was shown by Lloyd \cite{LloydPearson2002} in 2002.	
\end{remark}

\begin{pro}\label{pro13cl1} There exist parameters perturbations such that the following center family  
\begin{equation}
\begin{array}{lll}\label{komquar_c2}
\begin{aligned}
\dot{x}&=& x\,(y - 1)\,(4\,x - 8 + y)(2\,b\,x - a\,y + a - 2\,b - 3)/18,\\
\dot{y}&=& y\,(x - 1)\,(4\,x - 2 + y)(2\,b\,x - a\,y + a - 2\,b - 3)2/9,
\end{aligned}
\end{array}
\end{equation}	
can produce, at least, thirteen small limit cycles bifurcation from $(1,1).$      
\end{pro}
\begin{proof}
Note that the system \eqref{komquar_c2} is the system \eqref{komcub_c1} with $b_2=2$, and $a_{3}=1$ translated to the origin and multiplied by a straight line $-ay - bx + 1$. Then the system became 
	\begin{equation}
\begin{array}{lll}\label{komquar_c2_transladado}
\begin{aligned}
\dot{x}=& -(y-2\,x + 1)(1+3\,x/2)y\,(1-a\,y - b\,x),\\
\dot{y}=& \phantom{-} \, \, (2\,x - y + 1)(1-3\,y)\,x\,(1-a\,y - b\,x).
\end{aligned}
\end{array}
\end{equation}	
The system perturbated is given by 
\begin{equation*}
\begin{aligned}
\dot{x}=&-( y-2\,x + 1)(1+3\,x/2)y\,(1-a\,y - b\,x)+(1+3\,x/2)(a_{3, 0}\,x^3 + a_{2, 1}\,x^2\,y \\  \phantom{=}&+ a_{1, 2}\,x\,y^2 + a_{0, 3}\,y^3 + a_{2, 0}\,x^2 + a_{1, 1}\,x\,y + a_{0, 2}\,y^2),\\
\dot{y}=& \, \, \, \phantom{-} (2\,x - y + 1)(1-3\,y)\,x\,(1-a\,y - b\,x) + (1-3y)(b_{3, 0}\,x^3 +b_{2, 1}\,x^2\,y \\ \phantom{=}& + b_{1, 2}\,x\,y^2 + b_{0, 3}\,y^3 + b_{2, 0}\,x^2 + b_{1, 1}\,x\,y + b_{0, 2}\,y^2).
\end{aligned}
\end{equation*}
Computing the linear terms of the first Lyapunov Constants concerning perturbative parameters $\Lambda=\left\{a_{30}, a_{21}, a_{12}, a_{03}, a_{20}, a_{11}, a_{02}, b_{30}, b_{2,1}, b_{1,2}, b_{03},b_{20}, b_{11}, b_{02} \right\}$, we obtain that the rank of Jacobian Matrix in relation to $\Lambda$ is eleven. Due to size, we will show only the first two Lyapunov Constants.
\begin{equation*}
	\begin{aligned}
		L_1= &-\frac{1}{3}(aa_{1,1} + 3ab_{0,2} + ab_{2,0} + a_{0,2}b + 3a_{2,0}b + bb_{1,1}+ 2a_{0,2}- a_{1,1}+ a_{1,2}+ a_{2,0}+ 3a_{3,0}\\& - 7b_{0,2} + 3b_{0,3}- 2b_{1,1}+ b_{2,0}+ b_{2,1}), \\
		L_2= & -\frac{1}{90}(12a^3a_{1,1} + 72a^3b02 + 12a^3b_{2,0} + 48a^2a_{0,2}b + 36a^2a_{2,0}b + 48a^2bb_{1,1} + 24aa_{1,1}b^2 \\&  - 36ab^2b02 + 24ab^2b_{2,0} - 12a_{0,2}b^3 - 132b^4 + 24a^2a_{0,2} - 66a^2a_{1,1} + 12a^2a12 \\& - 450a^2b_{0,2} + 72a^2b_{0,3} - 24a^2b_{1,1} + 12a^2b_{2,1} - 150aa_{0,2}b + 36aa_{0,3}b + 172aa_{1,1}b \\& - 18aa_{2,0}b + 36aa_{2,1}b + 408abb_{0,2} - 78abb_{1,1} + 36abb_{1,2} + 28abb_{2,0} + 36abb_{3,0}  \\& + 106a_{0,2}b^2 + 12a_{1,1}b^2 - 12a_{1,2}b^2 + 330a_{2,0}b^2 + 168b^2b_{0,2} - 72b^2b_{0,3} + 154b^2b_{1,1}  \\& - 12b^2b_{2,0} - 12b^2b_{2,1} - 42a^2 - 180aa_{0,2} + 36aa_{0,3} - 458aa_{1,1} - 72aa_{1,2} - 36aa_{2,0}  \\& + 6aa_{2,1} - 810ab_{0,2} - 492ab_{0,3} + 36ab_{1,1} - 36ab_{1,2} - 446ab_{2,0} - 36ab_{2,1} - 66ab_{3,0} \\& - 132a_{0,2}b - 78a_{0,3}b - 202a_{1,1}b + 166a12b - 1499a_{2,0}b - 18a_{2,1}b + 360a_{3,0}b  \\& - 1072bb_{0,2} + 444bb_{0,3} - 772bb_{1,1} - 42bb_{1,2} + 58bb_{2,0} + 94bb_{2,1} + 18bb_{3,0} - 664a_{0,2} \\&  - 192a_{0,3} + 464a_{1,1} - 452a_{1,2} - 494a_{2,0} - 42a_{2,1} - 1659a_{3,0} + 3620b_{0,2} - 768b_{0,3} \\ & + 1000b_{1,1} + 48b_{1,2} - 62b_{2,0} - 452b_{2,1} - 102b_{3,0}).
    \end{aligned}		
\end{equation*}

We will consider $b_{20}=0, b_{21}=0,$ and  $b_{30}=0$. With the parameters, $a_{02},a_{03},a_{11},a_{12},a_{20}$, $a_{21},a_{30},b_{02},b_{03},b_{11},b_{12}$ the linear terms of the thirteen Lyapunov Constants have rank eleven. Making a linear change of coordinates we can write $L_{i}=u_{i}$, $1 \leq u \leq 10$ and 
\begin{equation*}
	\begin{array}{lll}
		L_{11} = f_{1}(a,b)u_{11}, \; \; \; \;
		L_{12} = f_{2}(a,b)u_{11},  \; \; \; \;
		L_{13} = f_{3}(a,b)u_{11},
	\end{array}
\end{equation*}
where $f_{1},f_{2},$ and $f_{3}$ are polynomials in which the numerator has  degree $64$, $65$, and $66$. The total number of monomials are, respectively, $171568, 190036,$ and $209989$. We can find a numerical solution $(a^{*}, b^{*})$ in the statement for the algebraic system such that $f_{1}(a^{*}, b^{*})=0, f_{2}(a^{*}, b^{*})=0$, and $f_{3}(a^{*}, b^{*})\neq0$. Moreover, the intersection is transversal due to the determinant of the Jacobian matrix at the intersection point  $(a^{*}, b^{*})$ to be $-3.11537939\times 10^{140}$.

To obtain an analytic proof we will use a Computer Assisted Proof with the help of the Theorem of Poincar\'e--Miranda and the Theorem of Circles of Gershgorin, for more details see \cite{daCruzTorreNovas2019}. We will use the Theorem of Poincar\'e--Miranda for the existence of the intersection point of $f_{1}$ and $f_{2}$, and we will use the Theorem of Circles of Gershgorin to prove the transversality. The technical both theorems also are used to check that at point $f_{3}$ is nonvanishing. We will fix a square $Q=[-h,h]$ with $h=10^{-10}$ and we do a rational affine change of coordinates such that a good rational approximation of $(a^{*}, b^{*})$ be inside $Q$. Then we have 
\begin{equation*}
	\begin{aligned}
	f_{1}(S_{0}^{-}) & \subset [-2.0000000158\times 10^{-10}, -1.999999993284\times10^{-10}],\\
    f_{1}(S_{0}^{+}) & \subset [2.00000000671\times10^{-10}, 1.999999984198\times10^{-10}],  \\ 
    f_{2}(S_{0}^{-}) & \subset [-2.0000000148\times10^{-10}, -2.00000000506\times10^{-10}], \\
    f_{2}(S_{0}^{+}) & \subset [1.9999999949\times10^{-10}, 1.9999999851\times10^{-10}], \\
    f_{3}(Q)         & \subset [0.99999999994,1.00000000005 ],
	\end{aligned}
\end{equation*}
and we have proved the existence of $(a^{*}, b^{*})$ such that $f_{3}$ is nonvanishing. In the computations, we have worked with rational numbers around $1400$ digits. To simplify the computations we have worked with the functions $f_{j}(a,b)=f_{j}(a,b)/f_{j}(0,0)$.

Finally, to prove this we need to check the transversality. Instead of computing the determinant of the Jacobian matrix of $(f_{1}, f_{2})$ concerning $(a^{*}, b^{*})$, we use the technical theorems already mentioned to get that, the elements in the Jacobian matrix for the transformed variables are, varying in $Q$, $J[1,1], J[2,2] \subset (0.999999,1.00000)$, $J[1,2], J[1,2] \subset (-1.30042\times10^{-19}, 1.20062\times10^{-18})$. Therefore, the determinant is different from zero. This finishes the proof.
\end{proof}

\begin{pro}\label{pro13lc2} There exist parameters perturbations such that the center given by 
	\begin{equation}
		\begin{array}{lll}\label{komquar_1}
			\begin{aligned}
				\dot{x}&=& \frac{1}{810}(29x^2 - 40xy + 40y^2 + 162x + 140y + 479)x(y - 1),\\
				\dot{y}&=&\frac{1}{810}(29x^2 - 40xy + 40y^2 + 162x + 140y + 479)y(x - 1),
			\end{aligned}
		\end{array}
	\end{equation}	
can produce, at least, thirteen small limit cycles bifurcation from $(1,1).$      
\end{pro}
\begin{proof}
\begin{figure}[h]
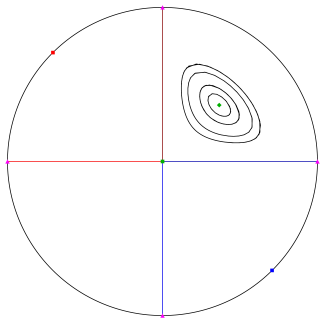
\vspace{0.3cm}
	\begin{overpic}[width=4.5cm]{centro_quartico}
	\end{overpic}
	\caption{ Phase portrait of the system~\eqref{komquar_1}.}
\end{figure}
 The system~(\ref{komquar_1}) has a center at the point $(1,1)$. We observe the system
  \begin{equation*}
 	\begin{array}{lll}\label{komquar_center1}
 		\begin{aligned}
 			\dot{x}=& \frac{x(y - 1)}{810}, \; \; \; \; \;
 			\dot{y}=&\frac{y(x - 1)}{810},
 		\end{aligned}
 	\end{array}
 \end{equation*}	
has a integrate factor $V=\frac{1}{xy}$, then the system has a center. As the system \eqref{komquar_1} is given by the system \ref{komquar_center1} multiplied by the quadratic curve of equilibrium points, follow that system \ref{komquar_1} has a center.  
  
Making a change of coordinates so that the origin is the center. Then, the unperturbed system is given by  
\begin{equation*}
	\begin{array}{lll}\label{komquar_2}
		\begin{aligned}
			\dot{x}=& -\left(\frac{29}{40}x^2 - xy + y^2 + x + y + 1\right)\left(\frac{9}{2}x + 1\right)y\\
			\dot{y}=& \phantom{-} \; \left(\frac{29}{40}x^2 - xy + y^2 + x + y + 1\right)\left(\frac{9}{2}y + 1\right)x,
		\end{aligned}
	\end{array}
\end{equation*}	
And the perturbated system is given by
\begin{equation*}
\begin{array}{lll}\label{komquar_3}
\dot{x} & = & -\left(\frac{29}{40}x^2 - xy + y^2 + x + y + 1\right)\left(\frac{9}{2}x + 1\right)y+(9x/2 + 1)(x^3a_{3, 0} + x^2ya_{2, 1}\vspace{0.5cm} \\ \vspace{0.5cm}
	 &\phantom{=}&+ xy^2a_{1, 2} + y^3a_{0, 3} + x^2a_{2, 0} + xya_{1, 1} + y^2a_{0, 2})\\
\dot{y} & = &  \left(\frac{29}{40}x^2 - xy + y^2 + x + y + 1\right)\left(\frac{9}{2}y + 1\right)x+(9y/2 + 1)(x^3b_{3, 0} + x^2yb_{2, 1}+ xy^2b_{1, 2} \vspace{0.5cm}\\ 
			&\phantom{=} &+ y^3b_{0, 3} + x^2b_{2, 0} + xyb_{1, 1} + y^2b_{0, 2}).
\end{array}
\end{equation*}	

To simplify the calculations we will take $b_{30}=0$, Calculating the linear terms of the first fourteen Lyapunov Constants, we obtain that the rank of Jacobian Matrix in relation with $\Lambda=\left\{a_{30}, a_{21},a_{12}, a_{03}, a_{20}, a_{11}, a_{02}, b_{30}, b_{2,1}, b_{1,2}, b_{03},b_{20}, b_{11}, b_{02} \right\}$ is eight, with the parameters $a_{0 2}, a_{0 3}, a_{1 1}, a_{1 2}, a_{2 0}, a_{2 1}, a_{3 0}, b_{03}$. We will show below the linear terms of the first three Lyapunov Constants.
\begin{equation*}
	\begin{aligned}
		L_1= & \phantom{-} \frac{1}{3}(-6a_{2,0} - 3a_{3,0} + a_{0,2} + a_{1,1} - a_{1,2} - 6b_{0,2} - 3b_{0,3} + b_{1,1} + b_{2,0} - b_{2,1}), \\
		L_2= & -\frac{1}{600}(6208a_{0,2} + 6109a_{1,1} - 6022a_{1,2} - 35970a_{2,0} - 14580a_{3,0} - 36597b_{0,2} \\& - 14712b_{0,3} + 6208b_{1,1} + 6109b_{2,0} - 6022b_{2,1}), \\
		L_{3} = &  \frac{1}{56000}(-23216498b_{2,1} + 24003392a_{0,2} - 15400a_{0,3} + 23569981a_{1,1} - 23216498a_{1,2}\\& - 138616730a_{2,0} - 55689120a_{3,0} - 141175473b_{0,2} - 56032276b_{0,3} + 24003392b_{1,1}\\& - 15400b_{1,2} + 23569981b_{2,0}).
     \end{aligned}		
\end{equation*}
We can write $L_{i}=u_{i}$ $1\leq u \leq 8,$ and make an appropriate change of coordinates so that we can eliminate the linear terms from $L_{9}$ until $L_{13}$, and we obtain five homogeneous equations of degree two in the parameters $b_{0 2}, b_{1 1}, b_{1 2}, b_{2 0},$ and $b_{2 1}$, when
\begin{equation*}
	\begin{aligned}
		L_9  = & \; b_{02} g_{1}(b_{02}, b_{11}, b_{12},b_{20}, b_{21}), \; \; \;
		L_{10} = & b_{02} g_{2}(b_{02}, b_{11}, b_{12},b_{20}, b_{21}),  \\
		L_{11} = & \; b_{02} g_{3}(b_{02}, b_{11}, b_{12},b_{20}, b_{21}), \; \; \;
		L_{12} = & b_{02} g_{4}(b_{02}, b_{11}, b_{12},b_{20}, b_{21}),  \\
		L_{13} = & \; b_{02} g_{5}(b_{02}, b_{11}, b_{12},b_{20}, b_{21}).   \\
	\end{aligned}		
\end{equation*}

Solving the equations $g_{i},$ $ 9\leq i \leq 12$, we find a point $p^{*}$ such that $L_{i}(p^{*})=0$ for $ 9\leq i \leq 12$ and we have $L_{13}(p^{*}) \neq 0$. Moreover, the interception among $g_{i},$ $ 9\leq i \leq 12$ at $p^{*}$ is transversal, that is, the Jacobian Matrix of $g_{i},$ $ 9\leq i \leq 12$ at $p^{*}$ is no zero. Therefore, there exist perturbative parameters such that thirteen small limit cycles bifurcating from the center.
\end{proof}

Finally, we will present two results for Kolmogorov quintic systems with the new lower bounds. The first, again, is based on the center used in Proposition~\ref{pro6cl}. And the second is based on the center used in Proposition~\ref{pro13lc2}, in which we find the best number of limit cycles for the quartic Kolmogorov system. We obtain these results as follows:

\begin{pro} For the center family of Kolmogorov quintic system,
\begin{equation}
\begin{array}{lll}\label{komquar_c3}
\begin{aligned}
\dot{x}&=& \frac{1}{768}(9x - 18 + y)x(y - 1)(9x^2 + 2y^2 - 18x - 4y - 21),\\
\dot{y}&=& \frac{1}{256}3(9x - 2 + y)y(x - 1)(9x^2 + 2y^2 - 18x - 4y - 21).
\end{aligned}
\end{array}
\end{equation}	
there exist polynomial perturbations of degree quintic such that at least 18 limit cycles bifurcation from $(1,1).$     
\end{pro}
\begin{proof}[Proof of Theorem~\ref{thmmainn_3}]
We observe that the system~(\ref{komquar_c3}) is the system~(\ref{komcub_c1}) multiplied by $(1-4x-2y)^2$, then system~(\ref{komquar_c3}) has a center at $(1,1)$. Calculating the linear terms of the first twenty-four Lyapunov Constants, we obtain that the rank of the Jacobian Matrix is eighteen in the parameters $a_{02},a_{03},a_{04},a_{11},a_{12},a_{13},a_{20},a_{21},a_{22},a_{30},a_{31},a_{40},b_{02},b_{03}$, $b_{04},b_{1,1},b_{12},b_{13}$. There is no more limit cycles using the terms of order two of Lyapunov Constants.
\end{proof}

\begin{pro}\label{pro22cl} For the center of Kolmogorov quintic system,
	\begin{equation}
		\begin{array}{lll}\label{komquint_c5}
			\begin{aligned}
				\dot{x}=& \;\phantom{-}\frac{1}{A}((-21 + 8x + 4y)(29x^2 - 40xy + 40y^2 + 162x + 140y + 479)x(y - 1)),\\
				\dot{y}=&-\frac{1}{A}((-21 + 8x + 4y)(29x^2 - 40xy + 40y^2 + 162x + 140y + 479)y(x - 1)).
			\end{aligned}
		\end{array}
	\end{equation}	
there exist polynomial perturbations of degree quintic such that at least 22 limit cycles bifurcation from $(1,1)$, where $A=7290$.   
\end{pro}

\begin{figure}[h]\vspace{0.3cm}
	\begin{overpic}[width=4cm]{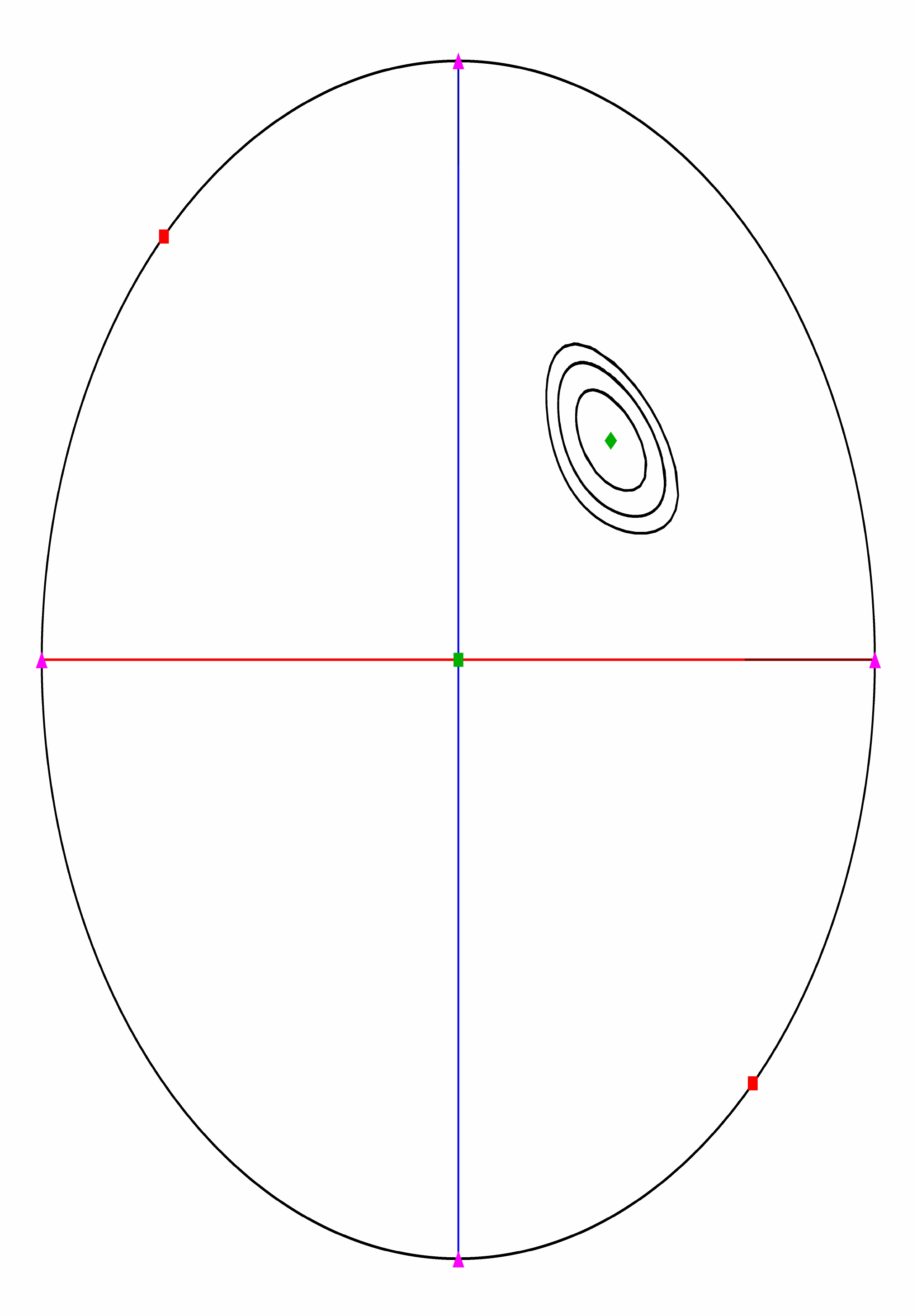}
	\end{overpic}
	\caption{ Phase portrait of the system~\eqref{komquint_c5}}
\end{figure}
\begin{proof}
Making a change of coordinates translating the center to origin, we obtain the following center
\begin{equation}
	\begin{array}{lll}\label{komquint_c6}
		\begin{aligned}
			\dot{x}&=&-(1 - 4x - 2y)(29/40x^2 - xy + y^2 + x + y + 1)((9x)/2 + 1)y,\\
			\dot{y}&=&(1 - 4x - 2y)(29/40x^2 - xy + y^2 + x + y + 1)((9y)/2 + 1)x.
		\end{aligned}
	\end{array}
\end{equation}	
We observe that the system \eqref{komquint_c6} is the system \eqref{komquar_c2_transladado} multiplied by $(1-4x-2y)^2$, then system \ref{komquint_c5} has a center at $(1,1)$. 

Calculating the linear terms of the first twenty-four Lyapunov Constants, we obtain that the rank of the Jacobian Matrix is thirteen in the parameters $a_{02},a_{03},a_{04},a_{11},a_{12},a_{13},a_{20}$, $a_{21},$ $a_{22},a_{30},a_{31},a_{40},b_{02}$. We will take the parameters $b_{04}=b_{11}=0$ to simplify the calculations. 

We can write $L_{i}=u_{i}$ $1\leq u \leq 13,$ and make an appropriate change of coordinates so that we can eliminate the linear terms from $L_{14}$ until $L_{22}$, and we obtain eight homogeneous equations of degree two. Solving the equations $L_{i},$ $ 14\leq i \leq 21$, we find a point $p^{*}$ such that $L_{i}(p^{*})=0$ for $ 14\leq i \leq 21$ and we have $L_{22}(p^{*}) \neq 0$. Moreover, the interception among $L_{i},$ $ 14\leq i \leq 22$ at $p^{*}$ is transversal, that is, the Jacobian Matrix of $L_{i},$ $ 14\leq i \leq 21$ at $p^{*}$ is no zero. Therefore, there exist perturbative parameters such that twenty-two small limit cycles bifurcating from the center.
\end{proof}

\section*{Acknowledgements}

Yagor Romano Carvalho was supported by São Paulo Paulo Research Foundation (FAPESP) grants number 2022/03800-7, 2021/14695-7. Leonardo P.C. da Cruz was supported by São Paulo Paulo Research Foundation (FAPESP) grants number 2022/14484-9, 2021/14987-8. Luiz F.S. Gouveia was supported by São Paulo Paulo Research Foundation (FAPESP) grants number 2022/03801-3, 2020/04717-0.

\section{Appendix}

\begin{teo}[\cite{gersh}, Circles of Gershgorin]
	Let $A=(a_{i,j}) \in \C^{n\times n}$ and $\alpha_{k}$ its eigenvalues. Consider for each $i=1,\ldots,n$ 
	\[D_{i}=\{ z \in \C : |z-a_{i,i}| \leq r_{i}\},\]
	where $r_{i}=\displaystyle\sum_{i\neq j}|a_{i,j}|$. So, for all $k$, each $\alpha_{k} \in D_{i}$ for some $i$.
\end{teo}

\begin{teo}[\cite{PM}, Poincar\'e--Miranda]
	Let $c$ be a positive real number and $S=[-c,c]^{n}$ a n-dimensional cube. Consider $f=(f_{1},\ldots,f_{n}):S \rightarrow \R^{n} $ a continuous function such $f_{i}(S_{i}^{-})<0$ and $f_{i}(S_{i}^{+})>0$ for each $i \leq n$, where $S_i^{\pm}=\left\{(x_{1},\ldots,x_{n}) \in S: x_{i}=\pm c \right\}$. So, there exists $d \in S$ such that $f(d)=0$. 
\end{teo}

%
%
%
\end{document}